\documentclass[12pt]{article}
\usepackage{theorem,amsmath,amssymb,latexsym}
\usepackage[all]{xy}

\usepackage{amscd}

\theoremstyle{plain}
\theoremheaderfont{\normalfont\bfseries}
\theorembodyfont{\itshape}

\newtheorem{theorem}{Theorem}[section]
\newtheorem{proposition}[theorem]{Proposition}
\newtheorem{lemma}[theorem]{Lemma}
\newtheorem{corollary}[theorem]{Corollary}

\theorembodyfont{\upshape}
{\theoremstyle{break}

}

\theorembodyfont{\upshape}

\theorembodyfont{\upshape}

\theorembodyfont{\upshape} 
\newtheorem{definition}[theorem]{Definition}

\theorembodyfont{\upshape}

\theorembodyfont{\itshape}

\theorembodyfont{\upshape}

\theorembodyfont{\upshape}

\theorembodyfont{\upshape}

\theorembodyfont{\upshape}

\theorembodyfont{\upshape}

\theorembodyfont{\upshape}
\newtheorem{proof}{Proof}

\theorembodyfont{\upshape}

\newtheorem{acknowledgement}{Acknowledgement}

\newcommand{\qed}{\nopagebreak\par\hspace*{\fill}$\square$\par\vskip2mm}
\newcommand{\myemail}[1]{\indent \emph{E-mail:} {\tt #1}}
\newcommand{\myaddress}[1]{\indent {\sc #1}\par}

\title{A characterization of ramification groups\\
of local fields with imperfect residue fields}
\author{Takeshi Saito}
\date{}

\begin{document}
\maketitle

\begin{abstract}
We give a characterization of ramification groups
of local fields with imperfect residue fields,
using those for local fields
with perfect residue fields.
As an application,
we reprove an equality
of ramification groups
for abelian extensions
defined in different ways.
\end{abstract}

Let $K$ be a henselian discrete
valuation field.
Let $L$ be a Galois extension
and let $G={\rm Gal}(L/K)$
be the Galois group.
In the classical case where
the residue field $F$ of $K$
is perfect,
the lower numbering filtration
$(G_{i,{\rm cl}})$ indexed by integers $i\geqq 0$
is defined by
$G_{i,{\rm cl}}={\rm Ker}(G\to
{\rm Aut}({\cal O}_L/
{\mathfrak m}^{i+1}_L))$.
Further, the upper numbering
filtration $(G^r_{\rm cl})$
indexed by rational numbers
$r>0$ is defined 
using the Herbrand function
to renumber the lower numbering
\cite[Chapitre IV, Section 3]{CL}.

In the general case where
the residue field need not be perfect,
an upper numbering filtration $(G^r)$
indexed by rational numbers
$r>0$ is defined first in \cite{AS}
using rigid geometry
and later in \cite{TJM}
purely in the language of schemes.
In the classical case where
$F$ is perfect,
they are related to each other
by the relation
$G^r=G^{r-1}_{\rm cl}$.

We give an axiomatic characterization
of the filtration $(G^r)$.
The axiom has two conditions.
The first condition is the relation in the classical case
above.
The second condition is the compatibility
with tangentially dominant extensions.
A similar approach
reducing to the classical case
was proposed
in \cite{Bo}.

For a discrete valuation field $K$,
the tangent space at an algebraic
closure $\bar F$ of the residue field $F$
is defined as an $\bar F$-vector space
using the cotangent complex. 
In the classical case where
the residue field is perfect,
it is nothing but the scalar
extension of the Zariski tangent
space defined as the dual
${\rm Hom}_F({\mathfrak m}_K/
{\mathfrak m}_K^2,\bar F)$.
An extension of discrete valuation
fields is said to be tangentially dominant
if the induced morphism
on the tangent spaces is dominant
(Definition \ref{dfcot}.2).
An unramified extension is 
tangentially dominant
and a tangentially dominant
extension has 
ramification index 1.

The uniqueness is a consequence of
the existence of tangentially dominant
extension with perfect residue field.
The existence follows from
the functorial properties of the filtration
$(G^r)$.

For $r>1$,
the graded quotient
${\rm Gr}^rG=G^r/G^{r+}$ is
defined by 
$G^{r+}=\bigcup_{s>r}G^s$
and is
an ${\mathbf F}_p$-vector space.
A canonical injection
\begin{equation}
{\rm Hom}({\rm Gr}^rG,{\mathbf F}_p)
\to
{\rm Hom}_{\bar F}(
{\mathfrak m}_{\bar K}^r/
{\mathfrak m}_{\bar K}^{r+},
H_1(L_{\bar F/S}))
\label{eq1}
\end{equation}
is defined in 
\cite[(4.20)]{red},
as a generalization of
a non-logarithmic variant of
the refined Swan conductor
defined by Kato in \cite{K}.
We also give
an axiomatic characterization of this
morphism, similar to the characterization
for $G^r$ itself.

As an application of the characterizations,
we give a new proof of the equality
of two filtrations for abelian extensions
in positive characteristic.
By the Hasse-Arf theorem,
the filtration $(G^n)$ defined in \cite{AS}
is in fact indexed by integers $n>1$
for abelian extensions.
The other filtration is the filtration 
$(G^n_{\rm Ma})$ defined
by Matsuda in \cite{M}
as a modification of that defined
by Kato in \cite{K}.
The equality was proved in \cite{aml}
except for the smallest index $n=2$
and the remaining case was
proved by Yatagawa in \cite{Y}.
The equality is proved by verifying
that the filtration $(G^n_{\rm Ma})$ 
satisfies the same axiom.
We also prove that the injection
(\ref{eq1})
equals the morphism
${\rm rsw}'$ defined in
\cite{M} and \cite{Y},
as a modification of
the refined Swan conductor
defined in \cite{K}.

A variant $(G^r_{\log})$
of the upper numbering filtration $(G^r)$
called the logarithmic 
upper numbering filtration
is also defined in \cite{AS}.
In the case where the ramification
index $e_{L/K}$ is 1,
the two filtrations are the same:
$G^r=G^r_{\log}$.
If $K'$ is a log smooth extension of
$K$ and $L'=LK'$,
the canonical injection
$G'={\rm Gal}(L'/K')\to 
G={\rm Gal}(L/K)$ is known to induce
isomorphisms
$G'^{er}_{\log}=G^r_{\log}$
for $e=e_{K'/K}$.
Further if the ramification
index $e_{K'/K}$ is divisible
by $e_{L/K}$, we have $e_{L'/K'}=1$.
Thus, a characterization of $(G^r)$
gives an indirect characterization
$(G^r_{\log})$.

\begin{acknowledgement}
The author thanks the referee for
careful reading and the suggestion
to include comments on the logarithmic filtration.
The research is partially supported
by Grant-in-Aid (B) 19H01780.
\end{acknowledgement}

\section{Totally ramified case}

Let $K$ be a henselian discrete
valuation field.
Let $L$ be a totally ramified
Galois extension of $K$
and let $G={\rm Gal}(L/K)$
be the Galois group.
For a rational number $r>1$,
the upper ramification group
$G^r$ defined in \cite[Definition 3.4]{AS}
equals the subgroup
defined in \cite[Chapitre IV, Section 3]{CL}
denoted $G^{r-1}_{\rm cl}$,
by \cite[Proposition 3.7 (3)]{AS}.

Assume that $L$ is wildly ramified
and 
let $r>1$ be the largest rational
number such that 
the subgroup $G^r$
of the wild inertia subgroup
$P\subset G$
is non-trivial.
Let $E$ be the residue field
and $e=e_{L/K}$
be the ramification index.
We give a description
of the canonical injection
\begin{equation}
G^{r\vee}={\rm Hom}_{{\mathbf F}_p}
(G^r,{\mathbf F}_p)
\to
{\rm Hom}_E(
{\mathfrak m}_L^{e(r-1)}/
{\mathfrak m}_L^{e(r-1)+1}
,E)
\label{eqchclE}
\end{equation}
for the ${\mathbf F}_p$-vector space
$G^r$,
in the case where $L$ is totally ramified
over $K$.
The injection (\ref{eqchclE})
is a special case of (\ref{eq1}).

We begin with a description
of extensions of vector spaces
over a field of characteristic $p>0$
by ${\mathbf F}_p$-vector
spaces.

\begin{lemma}\label{lmapol}
Let $F$ be a field of characteristic $p>0$.

{\rm 1.}
Let $G\subset F$
be a finite subgroup of the additive group.
Then, the polynomial
\begin{equation}
a_1=\dfrac{\prod_{\sigma\in G}(X-\sigma)}
{\prod_{\sigma\in G,\,
\sigma\neq 0}(-\sigma)}
\label{eqa1}
\end{equation}
$\in F[X]$ is a unique additive 
separable polynomial
such that the coefficient
of degree $1$ is $1$ and that
the sequence
\begin{equation}
\begin{CD}
0@>>> G
@>>>
{\mathbf G}_a
@>{a_1}>>
{\mathbf G}_a
@>>>0
\end{CD}
\label{eqapol}
\end{equation}
is exact.

{\rm 2. (\cite[Proposition 2.1.6
(2)$\Rightarrow$(3)]{red})}
Let $E$ be an $F$-vector space of finite
dimension and
let $0\to G\to H\to E\to 0$
be an extension of $E$ by
an ${\mathbf F}_p$-vector space $G$
of finite dimension, as smooth group
schemes over $F$.
Define a morphism
\begin{equation}
[H]\colon
G^\vee=
{\rm Hom}_{{\mathbf F}_p}
(G,{\mathbf F}_p)
\to 
{\rm Ext}(E,{\mathbf F}_p)
=E^\vee
={\rm Hom}_F(E,F)
\label{eqH}
\end{equation}
by sending a character
$\chi\colon G\to {\mathbf F}_p$
to the linear form $f\colon E\to F$
such that there exists
a commutative diagram
$$\begin{CD}
0@>>> G@>>> H@>>> E@>>>0\\
@.@V{\chi}VV@VVV@VVfV@.\\
0@>>>{\mathbf F}_p
@>>> {\mathbf G}_a
@>{x^p-x}>> {\mathbf G}_a
@>>> 0.
\end{CD}$$
If $H$ is connected, then
the morphism
$[H]\colon
G^\vee
\to E^\vee$ is an injection.
\end{lemma}

\begin{proof}
1.
By \cite[Lemma 2.1.5]{red},
$a=\prod_{\sigma\in G}(X-\sigma)
\in F[X]$ is an additive
separable polynomial
such that (\ref{eqapol})
with $a_1$ replaced by $a$
is exact.
Since the coefficient in $a$
of degree $1$ is
$\prod_{\sigma\in G,\,
\sigma\neq 0}(-\sigma)$,
the assertion follows. 
\qed
\end{proof}

\medskip
Let $K$ be a henselian discrete
valuation field and
$L$ be a totally ramified 
Galois extension
of degree $e$
of Galois group $G$.
Let $\alpha\in L$
be a uniformizer
and let $E=F$ denote the residue field.
The minimal polynomial
$f\in {\cal O}_K[X]$
is an Eisenstein polynomial
and the constant term $\pi=f(0)$
is a uniformizer of $K$.
We define a closed immersion
$T={\rm Spec}\, {\cal O}_L
\to Q={\rm Spec}\, {\cal O}_K[X]$
by sending $X$ to $\alpha$.
For a rational number $r>1$ 
such that $er\in{\mathbf Z}$,
define a dilatation
$$Q^{[r]}_T=
{\rm Spec}\, {\cal O}_L[X]
\Bigl[\dfrac f{\alpha^{er}}\Bigr]
\to Q_T=
{\rm Spec}\, {\cal O}_L[X].$$
The generator
$f$ of the kernel
$I={\rm Ker}({\cal O}_K[X]
\to {\cal O}_L)$
defines a basis over ${\cal O}_L$
of the conormal module
$N_{T/Q}=I/I^2$
and $\alpha^{er}$
defines a basis
of the $E$-vector space
${\mathfrak m}_L^{er}/
{\mathfrak m}_L^{er+1}$.
As subspaces of
$N_{E/Q}=J/J^2$
for $J={\rm Ker}({\cal O}_K[X]
\to E)=
(X, f)
=(X,\pi)$,
we have an equality
\begin{equation}
N_{T/Q}\otimes_{{\cal O}_L}E
={\mathfrak m}_K/
{\mathfrak m}_K^2
\label{eqNTQ}
\end{equation}
since $f$ is an Eisenstein polynomial.
The basis $f$ of 
$N_{T/Q}\otimes_{{\cal O}_L}E$
corresponds to
the uniformizer
$\pi=f(0)
\in {\mathfrak m}_K/
{\mathfrak m}_K^2$.

By sending
$S$ to $f/\alpha^{er}$,
we define an isomorphism
${\cal O}_L[X,S]/
(f-\alpha^{er}S)
\to {\cal O}_L[X]
[f/\alpha^{er}]$.
Since $f$ is an Eisenstein polynomial,
the reduced closed fiber
$Q^{[r]}_E=
{\rm Spec}\, ({\cal O}_L[X]
[f/\alpha^{er}]
\otimes_{{\cal O}_L}E)_{\rm red}$
is identified with
${\rm Spec}\, E[S]$.
By this identification
and (\ref{eqNTQ}),
we define an isomorphism
\begin{align}
Q^{[r]}_E\to
{\rm Hom}_E({\mathfrak m}_L^{er}
/{\mathfrak m}_L^{er+1},
N_{T/Q}\otimes
_{{\cal O}_L}E)^\vee
\to&\, 
{\rm Hom}_E({\mathfrak m}_L^{er}
/{\mathfrak m}_L^{er+1},
{\mathfrak m}_K/
{\mathfrak m}_K^2)^\vee
\nonumber
\\
&
=
{\mathfrak m}_L^{e(r-1)}
/{\mathfrak m}_L^{e(r-1)+1}
\label{eqQr}
\end{align}
of smooth group schemes of
dimension 1 over $E$.

Let $Q^{(r)}_T\to
Q^{[r]}_T$ be the normalization
and define a section
$T\to Q^{(r)}_T$
to be the unique lifting of
the section
$T\to Q_T$
defined by sending $X$ to $\alpha$.
Let $Q^{(r)}_E$
denote the reduced part
of the closed fiber
$Q^{(r)}_T\times_T{\rm Spec}\, E$
and let
$Q^{(r)\circ}_E\subset Q^{(r)}_E$
denote the connected 
component containing
the image of the closed point
of $T$ by the section
$T\to Q^{(r)}_T$.

\begin{proposition}\label{prcl}
Let $K$ be a henselian discrete
valuation field with residue field
$F$ of characteristic $p>0$.
Let $L$ be a totally ramified 
Galois extension
of degree $n=e$
with residue field $E=F$
and let $G={\rm Gal}(L/K)$
be the Galois group.
Let $\alpha\in L$
be a uniformizer
and let $f\in {\cal O}_K[X]$
be the minimal polynomial.
Decompose
$f=\prod_{i=1}^n(X-\alpha_i)$
so that $\alpha_n=\alpha$
and ${\rm ord}_L(\alpha_i-\alpha_n)$
is increasing in $i$.

{\rm 1.}
Let $r>1$ be the largest rational number
such that $G^r\neq 1$.
Then, we have 
\begin{equation}
er=
{\rm ord}_Lf'(\alpha)
+{\rm ord}_L(\alpha_{n-1}-\alpha_n).
\label{eqr}
\end{equation}
Define an injection
$\beta\colon G^r\to {\mathbf G}_a$
by $\beta(\sigma)
\equiv\dfrac{\sigma(\alpha)-\alpha}
{\alpha_{n-1}-\alpha_n}
\bmod {\mathfrak m}_L$
and
an additive polynomial
$b_1\in E[X]$
by $b_1=
\prod_{\sigma\in G^r}
(X-\beta(\sigma))
/\prod_{\sigma\in G^r,\sigma\neq 1}
(-\beta(\sigma))$.
Define an isomorphism
${\mathbf G}_a\to
{\mathfrak m}_L^{e(r-1)}
/{\mathfrak m}_L^{e(r-1)+1}$ by
$f'(\alpha)
(\alpha_{n-1}-\alpha_n)/f(0)
\in {\mathfrak m}_L^{e(r-1)}$
and identify
$Q^{[r]}_E$ with
${\mathfrak m}_L^{e(r-1)}
/{\mathfrak m}_L^{e(r-1)+1}$
by the isomorphism
{\rm (\ref{eqQr})}.
Then, there exists an isomorphism
\begin{equation}
\begin{CD}
0@>>>
G^r@>>>
Q^{(r)\circ}_E
@>>>
Q^{[r]}_E
@>>>
0
\\
@.@|@AAA@AA{f'(\alpha)
(\alpha_{n-1}-\alpha_n)/f(0)}A@.\\
0@>>>
G^r@>\beta>>
{\mathbf G}_a@>{b_1}>>
{\mathbf G}_a@>>>0
\end{CD}
\label{eqGrcl}
\end{equation}
of exact sequences.

{\rm 2.}
Let $i>0$ be
the largest integer such that
$G_{i,{\rm cl}}={\rm Ker}(G\to
{\rm Aut}({\cal O}_L/{\mathfrak m}_L^{i+1}))
\neq 1$.
Then, we have
\begin{equation}
i=
{\rm ord}_L(\alpha_{n-1}-\alpha_n)-1.
\label{eqi}
\end{equation}

Let $K\subset M\subset L$
be the intermediate extension
corresponding to $G_{i,{\rm cl}}\subset G$
and let $U^i_L=1+{\mathfrak m}^i_L
\subset L^\times$
and $U^i_M=1+{\mathfrak m}^i_M
\subset M^\times$ be the 
multiplicative subgroups.
Let
$N^i\colon U^i_L/U^{i+1}_L
\to 
U^i_M/U^{i+1}_M$
denote the morphism induced by the norm
$N_{L/M}\colon L^\times \to M^\times$
and
$T^i\colon U^i_L/U^{i+1}_L
={\mathfrak m}^i_L/{\mathfrak m}^{i+1}_L
\to
U^i_M/U^{i+1}_M
={\mathfrak m}^i_M/{\mathfrak m}^{i+1}_M$
be the isomorphism
induced by the trace
${\rm Tr}_{L/M}\colon L\to M$.
Define an isomorphism
${\mathbf G}_a\to
U^i_L/U^{i+1}_L$
by sending $1$
to the class of
$\alpha_{n-1}/\alpha_n\in U^i_L$.
Then, the diagram
\begin{equation}
\begin{CD}
0@>>>
G_{i,{\rm cl}}@>{\sigma
\mapsto \sigma(\alpha)/\alpha}>>
U^i_L/U^{i+1}_L
@>{(T^i)^{-1}\circ N^i}>>
U^i_L/U^{i+1}_L
@>>>
0
\\
@.@|@AA{\alpha_{n-1}/\alpha_n}A
@AA{\alpha_{n-1}/\alpha_n}A@.\\
0@>>>
G_{i,{\rm cl}}@>\beta>>
{\mathbf G}_a@>{b_1}>>
{\mathbf G}_a@>>>0
\end{CD}
\label{eqGicl}
\end{equation}
is an isomorphism of exact sequences.
\end{proposition}

\begin{proof}
1.
We have (\ref{eqr}) by
\cite[Lemma 3.3.1.5]{red}.

We have a commutative diagram
(\ref{eqGrcl}) with $b_1$ and
$f'(\alpha)
(\alpha_{n-1}-\alpha_n)$
replaced by
$b=\prod_{\sigma\in G^r}
(X-\beta(\sigma))$
and $c=
\prod_{i=1}^m
(\alpha_n-\alpha_i)\cdot
(\alpha_{n-1}-\alpha_n)^{n-m}$
for $m=\#G-\#G^r$
by \cite[Lemma 3.3.1.1]{red},
since the canonical isomorphism
$N_{T/Q}\to N_{E/S}\otimes E$
maps $f$ to $f(0)$.
Since $b=
\prod_{\sigma\in G^r,\sigma\neq 1}
(-\beta(\sigma))\cdot b_1$
and
$c=\prod_{\sigma\in G^r,\sigma\neq 1}
(-\beta(\sigma))\cdot f'(\alpha)
(\alpha_{n-1}-\alpha_n)$,
we obtain (\ref{eqGrcl}).

2.
Since ${\cal O}_L=
{\cal O}_K[\alpha]$
and ${\rm ord}_L(\alpha_i-\alpha_n)$
is increasing,
the equality (\ref{eqi}) follows from
the definition of $G_{i,{\rm cl}}$.

By \cite[Chapitre V, Proposition 8, Section 6]{CL},
the morphism
$(T^i)^{-1}\circ N^i\colon
U^i_L/U^{i+1}_L\to
U^i_L/U^{i+1}_L$
is defined by a separable additive
polynomial such that
the coefficient of degree 1 is 1
and the upper line of (\ref{eqGicl})
is exact.
Since $\sigma(\alpha)/\alpha
=1+(\sigma(\alpha)-\alpha_n)/
(\alpha_{n-1}-\alpha_n)\cdot(
\alpha_{n-1}/\alpha_n-1)$,
the left square is commutative.
Since the left square is commutative,
the right square is also commutative
by the uniqueness of $b_1$.
\qed
\end{proof}

\begin{corollary}\label{corcl}
{\rm 1.}
We have 
\begin{equation}
er=
{\rm ord}_Lf'(\alpha)+(i+1).
\label{eqri}
\end{equation}

{\rm 2.}
There exists an isomorphism
\begin{equation}
\begin{CD}
0@>>>
G^r@>>>
Q^{(r)\circ}_E
@>>>
{\mathfrak m}_L^{e(r-1)}
/{\mathfrak m}_L^{e(r-1)+1}
@>>>
0
\\
@.@|@AAA@AA{f'(\alpha)\cdot
\alpha_n/f(0)}A@.\\
0@>>>
G_{i,{\rm cl}}@>>>
U^i_L/U^{i+1}_L
@>{(T^i)^{-1}\circ N^i}>>
U^i_L/U^{i+1}_L
@>>>
0
\end{CD}
\label{eqGricl}
\end{equation}
of exact sequences.
\end{corollary}

\begin{proof}
1.
The equality (\ref{eqri})
follows from (\ref{eqr})
and (\ref{eqi}).

2.
Combining (\ref{eqGrcl})
and (\ref{eqGicl}),
we obtain (\ref{eqGricl}).
\qed
\end{proof}

\medskip

By Lemma \ref{lmapol}.2,
the extension
in the upper line of (\ref{eqGricl})
defines a canonical injection
\begin{equation}
G^{r\vee}={\rm Hom}_{{\mathbf F}_p}
(G^r,{\mathbf F}_p)
\to
{\rm Hom}_E(
{\mathfrak m}_L^{e(r-1)}/
{\mathfrak m}_L^{e(r-1)+1}
,
E).
\label{eqchcl}
\end{equation}

Assume that the residue
field of $F$ is perfect
and let $L$ be a Galois extension
of $K$.
Let $K^{\rm ur}\subset L$
denote the maximum
unramified extension
corresponding
to the inertia subgroup $I\subset G$.
For a rational number $r>1$,
we apply the construction
of (\ref{eqchcl})
to the totally ramified
extension
$M\subset L$ of $K^{\rm ur}$ 
corresponding to
$G^{r+}=\bigcup_{s>r}
G^s\subset I\subset G$
and to $H^r={\rm Gr}^rG=
G^r/G^{r+}\subset H
={\rm Gal}(M/K^{\rm ur})
=I/G^{r+}$.
Let $e'=e_{M/K}$ be the ramification
index and $E'\subset E$ be
the residue field of $M$.
We obtain an injection
\begin{align}
({\rm Gr}^rG)^\vee=
{\rm Hom}_{{\mathbf F}_p}
({\rm Gr}^rG,{\mathbf F}_p)
\to&\,
{\rm Hom}_{E'}(
{\mathfrak m}_M^{e'(r-1)}/
{\mathfrak m}_M^{e'(r-1)+1},E')
\label{eqgrch}
\\
&
\subset
{\rm Hom}_E(
{\mathfrak m}_L^{e(r-1)}/
{\mathfrak m}_L^{e(r-1)+1},E).
\nonumber
\end{align}

For abelian extensions,
we have the Hasse-Arf theorem.

\begin{theorem}{\rm ({\cite[Chapitre V,
Section 7, Th\'eor\`eme 1]{CL}})}
\label{thmHA}

Let $K$ be a henselian discrete valuation
field with perfect residue field
and let $L$ be a finite abelian
extension of $K$.
Let $n\geqq 1$ be an integer
and $r$ be a rational number 
satisfying $n<r\leqq n+1$.
Then, we have
$G^r=G^{n+1}$.
\end{theorem}

\section{Tangent spaces
and a characterization 
of ramification groups}

\begin{definition}{\rm({\cite[Definition 1.1.8]{red}})}
\label{dfcot}
Let $K$ be a discrete valuation 
field,
$S={\rm Spec}\, {\cal O}_K$
and $F$ be the residue field.

{\rm 1.}
For an extension $E$ of $F$,
let $L_{E/S}$ denote the cotangent complex
and 
we call the spectrum
\begin{equation}
\Theta_{K,E}=
{\rm Spec}\, S(H_1(L_{E/S}))
\label{eqThA}
\end{equation}
of the symmetric algebra
over $E$ 
the tangent space of $S$ at $E$.

{\rm 2.}
If ${\cal O}_K\to {\cal O}_{K'}$
is a faithfully flat morphism of
discrete valuation rings,
we say that $K'$ is
an extension of
discrete valuation fields of $K$.
We say that
an extension $K'$
of 
discrete valuation fields of $K$
is tangentially dominant if,
for a morphism $\bar F\to \bar F'$
of algebraic closures of the residue fields,
the morphism
$$
S(H_1(L_{\bar F/S}))
\to
S(H_1(L_{\bar F'/S'}))$$
is an injection.
\end{definition}

The morphism
\begin{equation}
{\mathfrak m}_K/
{\mathfrak m}_K^2\otimes_F\bar F
=
H_1(L_{F/S})\otimes_F\bar F
\to 
H_1(L_{\bar F/S})
\label{eqLFm}
\end{equation}
defined by the functoriality
of cotangent complexes
is an injection
by \cite[Proposition 1.1.3.1]{red}.
The injection (\ref{eqLFm})
is an isomorphism if
$F$ is perfect.
The distinguished triangle
$L_{S/{\mathbf Z}}
\otimes^L_{{\cal O}_S}
\bar F
\to L_{\bar F/{\mathbf Z}}
\to L_{\bar F/S}\to $
defines a canonical surjection
\begin{equation}
H_1(L_{\bar F/S})
\to
\Omega^1_{{\cal O}_K}
\otimes_{{\cal O}_K}
\bar F
\label{eqLHO}
\end{equation}
\cite[Proposition 1.1.7.3]{red}
such that the composition with
(\ref{eqLFm}) is induced by
$d\colon 
{\mathfrak m}_K/
{\mathfrak m}_K^2\to
\Omega^1_{{\cal O}_K}
\otimes_{{\cal O}_K}F$.
If $K$ is of characteristic $p>0$,
(\ref{eqLHO})
is an isomorphism
by 
\cite[Proposition 1.1.7.3]{red}.
If $K'$
is a tangentially dominant
extension of $K$,
the morphism
$H_1(L_{\bar F/S})
\to
H_1(L_{\bar F'/S'})$
is an injection.

\begin{proposition}{\rm ({\cite[Proposition 1.1.10]{red}})}
\label{prcot}
Let $K\to K'$
be an extension of discrete valuation
fields.
We consider the following conditions:

{\rm (1)}
The ramification index 
$e_{K'/K}$ is $1$
and $F'={\cal O}_{K'}/
{\mathfrak m}_{K'}$
is a separable extension of 
$F={\cal O}_K/
{\mathfrak m}_K$.

{\rm (2)}
The extension $K'$
is tangentially dominant over $K$.

{\rm (3)}
The ramification index 
$e_{K'/K}$ is $1$.

\noindent
Then, we have
the implications
{\rm (1)}$\Rightarrow${\rm (2)}$\Rightarrow${\rm (3)}.
\end{proposition}

\begin{theorem}\label{thmfil}
Let $r>1$ be a rational number.
For finite Galois extensions $L$
of henselian discrete valuation fields $K$,
there exists a unique way
to define a normal subgroup
$G^r$ of the Galois group
$G={\rm Gal}(L/K)$
satisfying the following conditions:

{\rm (1)}
If the residue field of $K$
is perfect, then
$G^r=G^{r-1}_{\rm cl}$.

{\rm (2)}
Let $K'$ be a tangentially dominant
extension of $K$.
Then the natural injection
$G'={\rm Gal}(L'/K')
\to G$ for $L'=LK'$
induces an isomorphism
$G^{\prime r}\to G^r$.
\end{theorem}

For a separable closure
$\bar K$ of $K$,
extend the normalized discrete
valuation ${\rm ord}_K$
to $\bar K$.
For a rational number $r$,
set
${\mathfrak m}_{\bar K}^r
=\{x\in \bar K\mid {\rm ord}_K
x\geqq r\}
\supset
{\mathfrak m}_{\bar K}^{r+}
=\{x\in \bar K\mid {\rm ord}_K
x> r\}$.
The quotient
${\mathfrak m}_{\bar K}^r
/{\mathfrak m}_{\bar K}^{r+}$
is a vector space of dimension 1
over the residue field $\bar F$.
For $r>1$,
define $G^{r+}=\bigcup_{s>r}G^s$
and ${\rm Gr}^rG=G^r/G^{r+}$.

\begin{theorem}\label{thmchar}
Let $r>1$ be a rational number.
For finite Galois extensions $L$
of henselian discrete valuation fields $K$,
for morphisms
$L\to \bar K$
to separable closures over $K$
and for the residue field $\bar F$
of $\bar K$,
there exists a unique way
to define an injection
\begin{equation}
{\rm Hom}({\rm Gr}^rG,{\mathbf F}_p)
\to
{\rm Hom}_{\bar F}(
{\mathfrak m}_{\bar K}^r/
{\mathfrak m}_{\bar K}^{r+},
H_1(L_{\bar F/S})).
\label{eqgrL}
\end{equation}
satisfying the following conditions:

{\rm (1)}
Assume that the residue field of $K$
is perfect.
Let $E$ be the residue field of $L$,
$e=e_{L/K}$ be
the ramification index
and identify
$
{\rm Hom}_E(
{\mathfrak m}_L^{e(r-1)}/
{\mathfrak m}_L^{e(r-1)+1},
E)
$
with a subgroup of
$
{\rm Hom}_{\bar F}(
{\mathfrak m}_{\bar K}^r/
{\mathfrak m}_{\bar K}^{r+},
H_1(L_{\bar F/S}))$
by the injection
${\mathfrak m}_K/
{\mathfrak m}_K^2
\to
H_1(L_{\bar F/S})$
{\rm (\ref{eqLFm})}.
Then, the diagram
\begin{equation}
\begin{CD}
{\rm Hom}({\rm Gr}^rG,{\mathbf F}_p)
@>>>
{\rm Hom}_{\bar F}(
{\mathfrak m}_{\bar K}^r/
{\mathfrak m}_{\bar K}^{r+},
H_1(L_{\bar F/S}))\\
@|@AAA\\
{\rm Hom}({\rm Gr}^rG,{\mathbf F}_p)
@>{\rm (\ref{eqchcl})}>>
{\rm Hom}_E(
{\mathfrak m}_L^{e(r-1)}/
{\mathfrak m}_L^{e(r-1)+1},
E)
\end{CD}
\label{eqgrcl}
\end{equation}
is commutative.

{\rm (2)}
Let $K'$ be a tangentially dominant
extension of $K$,
let $\bar K\to \bar K'$
be a morphism of separable
closures extending
$L\to L'=LK'$
and let $\bar F\to \bar F'$
be the morphism of residue fields.
Then, for the natural injection
$G'={\rm Gal}(L'/K')
\to G$,
the diagram
\begin{equation}
\begin{CD}
{\rm Hom}({\rm Gr}^rG,{\mathbf F}_p)
@>>>
{\rm Hom}_{\bar F}(
{\mathfrak m}_{\bar K}^r/
{\mathfrak m}_{\bar K}^{r+},
H_1(L_{\bar F/S}))\\
@VVV@VVV\\
{\rm Hom}({\rm Gr}^rG',{\mathbf F}_p)
@>>>
{\rm Hom}_{\bar F'}(
{\mathfrak m}_{\bar K'}^r/
{\mathfrak m}_{\bar K'}^{r+},
H_1(L_{\bar F'/S'}))
\end{CD}
\label{eqLL}
\end{equation}
is commutative.
\end{theorem}

The uniqueness is a consequence
of the following
existence of a tangentially dominant
extension with perfect residue field.

\begin{proposition}{\rm (\cite[Proposition 1.1.12]{red})}\label{prperf}
Let $K$ be a discrete valuation field.
Then, there exists a
tangentially dominant extension
$K'$ of $K$ such that
the residue field $F'$ is perfect.
\end{proposition}

\begin{proof}\hspace{-2mm}{\bf of Theorem \ref{thmfil}}
We show the uniqueness.
By Proposition \ref{prperf},
there exists a tangentially dominant
extension $K'$
of $K$ with perfect residue field.
Let $G'={\rm Gal}(L'/K')
\to G$ be the natural injection
for $L'=LK'$.
Then, by the conditions (1) and (2),
the subgroup $G^r\subset G$
is the image of $G^{\prime r-1}_{\rm cl}
\subset G'$.

To show the existence,
it suffices to prove that the subgroup
$G^r\subset G$ defined in \cite{AS}
satisfies the conditions (1) and (2).
The equality $G^r=G^{r-1}_{\rm cl}$
is proved in
\cite[Proposition 3.7 (3)]{AS}.
The condition (2) is satisfied by
\cite[Proposition 4.2.4 (1)]{red}.
\qed
\end{proof}

\begin{proof} \hspace{-2mm}{\bf of Theorem \ref{thmchar}}
We show the uniqueness.
If the residue field is perfect,
the morphism (\ref{eqgrL})
is uniquely determined by the commutative
diagram (\ref{eqgrcl})
since its right vertical arrow is an injection 
induced by the injection (\ref{eqLFm}).
In general,
by Proposition \ref{prperf},
there exists a tangentially dominant
extension $K'$
of $K$ with perfect residue field.
Then, the morphism (\ref{eqgrL})
is uniquely determined by the commutative
diagram (\ref{eqLL})
since its right vertical arrow 
is an injection.

To show the existence,
it suffices to prove that
the morphism
\cite[(4.20)]{red} satisfies
the conditions (1) and (2). 
Assume that 
the residue field is perfect.
To show the commutative diagram
(\ref{eqgrcl}),
we may assume that
$G^{r+}=1$
and ${\rm Gr}^rG=G^r$
by the construction of the morphisms.
Then, since
the construction of
{\rm (\ref{eqchcl})}
is a special case of
\cite[(4.20)]{red},
the condition (1) is satisfied.
The condition (2) follows from
\cite[(4.19)]{red}.
\qed
\end{proof}

\section{Abelian extensions}

\begin{theorem}\label{thmfila}
Let $r>1$ be a rational number.

{\rm 1.}
For finite abelian extensions $L$
of henselian discrete valuation fields $K$,
there exists a unique way
to define a normal subgroup
$G^r$ of the Galois group
$G={\rm Gal}(L/K)$
satisfying the following conditions:

{\rm (1)}
If the residue field of $K$
is perfect, then
$G^r=G^{r-1}_{\rm cl}$.

{\rm (2)}
Let $K'$ be a tangentially dominant
extension of $K$.
Then the natural injection
$G'={\rm Gal}(L'/K')
\to G$ for $L'=LK'$
induces an isomorphism
$G^{\prime r}\to G^r$.

{\rm 2.}
Let $L$ be a finite abelian extension
of a henselian discrete valuation field $K$
and let $n\geqq 1$ be the integer
satisfying $n<r\leqq n+1$.
Then, we have $G^r=G^{n+1}$.
\end{theorem}

\begin{theorem}\label{thmchara}
Let $n>1$ be an integer.
For finite abelian extensions $L$
of henselian discrete valuation fields $K$,
for morphisms
$L\to \bar K$
to separable closures over $K$
and for the residue fields $\bar F$ of 
$\bar K$,
there exists a unique way
to define an injection
\begin{equation}
{\rm Hom}(G^n/G^{n+1},{\mathbf F}_p)
\to
{\rm Hom}_F(
{\mathfrak m}_K^n/
{\mathfrak m}_K^{n+1},
H_1(L_{\bar F/S})).
\label{eqgrLa}
\end{equation}
satisfying the following conditions:

{\rm (1)}
Assume that the residue field of $K$
is perfect
and let the notation be as
in Theorem {\rm \ref{thmchar} (1)}.
Then the diagram
\begin{equation}
\begin{CD}
{\rm Hom}(G^n/G^{n+1},{\mathbf F}_p)
@>>>
{\rm Hom}_F(
{\mathfrak m}_K^n/
{\mathfrak m}_K^{n+1},
H_1(L_{\bar F/S}))\\
@|@AAA\\
{\rm Hom}(G^{n-1}_{\rm cl}
/G^n_{\rm cl},{\mathbf F}_p)
@>>>
{\rm Hom}_F(
{\mathfrak m}_K^{n-1}/
{\mathfrak m}_K^{n},
E)
\end{CD}
\label{eqgrcla}
\end{equation}
is commutative.

{\rm (2)}
Let $K'$ be a tangentially dominant
extension of $K$
and let the notation be as
in Theorem {\rm \ref{thmchar} (2)}.
Then, the diagram
\begin{equation}
\begin{CD}
{\rm Hom}(G^n/G^{n+1},{\mathbf F}_p)
@>>>
{\rm Hom}_F(
{\mathfrak m}_K^n/
{\mathfrak m}_K^{n+1},
H_1(L_{\bar F/S}))
\\
@VVV@VVV
\\
{\rm Hom}(G^{\prime n}/G^{\prime n+1},
{\mathbf F}_p)
@>>>
{\rm Hom}_{F'}(
{\mathfrak m}_{K'}^n/
{\mathfrak m}_{K'}^{n+1},
H_1(L_{\bar F'/S'}))
\end{CD}
\label{eqLLa}
\end{equation}
is commutative.
\end{theorem}

\begin{proof} \hspace{-2mm}{\bf of Theorem \ref{thmfila}}
1. is proved in the same way as
Theorem \ref{thmfil}.

2. By 1, this
follows from the Hasse-Arf theorem
Theorem \ref{thmHA}.
\qed
\end{proof}

\begin{proof} \hspace{-2mm}{\bf of Theorem \ref{thmchara}}
This is a special case of
Theorem \ref{thmchar}.
\qed
\end{proof}

\medskip

Assume that $K$ is
a henselian discrete valuation field
of equal characteristic $p>0$ and
let $L$ be a finite abelian extension.
Then, by the Hasse-Arf theorem
Theorem \ref{thmfila}.2
and by the isomorphism
$H_1(L_{\bar F/S})\to
\Omega^1_{{\cal O}_K}
\otimes_{{\cal O}_K}\bar F$
(\ref{eqLHO}),
for an integer $n>0$,
the injection 
{\rm (\ref{eqgrLa})}
defines an injection
\begin{equation}
{\rm Hom}(G^n/G^{n+1},{\mathbf F}_p)
\to
{\rm Hom}_F(
{\mathfrak m}_K^n/
{\mathfrak m}_K^{n+1},
\Omega^1_{{\cal O}_K}
\otimes_{{\cal O}_K}\bar F).
\label{eqgrM}
\end{equation}
A decreasing filtration $(G^n_{\rm Ma})$
indexed by integers $n>0$
is defined in \cite[Definition 3.1.1]{M}
as a non-logarithmic modification
of a filtration $(G^n_{\rm Ka})$
defined in \cite[Definition (2.1)]{K}.
Further, for an integer $n>0$,
a canonical morphism
\begin{equation}
{\rm rsw}'\colon
{\rm Hom}(G^n_{\rm Ma}/G_{\rm Ma}^{n+1},{\mathbf F}_p)
\to
{\rm Hom}_F(
{\mathfrak m}_K^n/
{\mathfrak m}_K^{n+1},
\Omega^1_{{\cal O}_K}
\otimes_{{\cal O}_K}\bar F)
\label{eqgrMY}
\end{equation}
is defined in 
{\rm \cite[Definition 3.2.5]{M}}
except the case
$p=2$, $n=2$
and in 
{\rm \cite[Definition 1.18]{Y}}
in the exceptional case
$p=2$, $n=2$,
as a modification of
the refined Swan conductor
defined in \cite[Corollary (5.2)]{K}.

As an application,
we give a new proof of
the equalities of the two filtrations
and the two morphisms,
different from that in
\cite{aml} and \cite{Y}.

\begin{corollary}\label{corab}
Let $L$ be an abelian extension
of a henselian discrete valuation field $K$
of equal characteristic $p>0$
and let $n>1$ be an integer.

{\rm 1.({\cite[Th\'eor\`eme 9.10 (i)]{aml}
for $n>1$,
\cite[Theorem 3.1]{Y}}
 for $n$ general)}
We have an equality
$G^n=G^n_{\rm Ma}$
of subgroups of $G$.

{\rm 2.(\cite[Th\'eor\`eme 9.10
(ii)]{aml} for $n>1$,
\cite[Corollary 2.13]{Y} for $n$ general)}
The injection
{\rm (\ref{eqgrM})}
is the same as  ${\rm rsw}'$
{\rm (\ref{eqgrMY})}.
\end{corollary}

The following proof
is by the reduction
to the logarithmic variant
\cite[Th\'eor\`eme 9.11]{aml} 
in the classical case
where the residue field
is perfect.

\begin{proof}
It suffices to show that the filtration 
$(G^n_{\rm Ma})$
and the morphism
${\rm rsw}'$
satisfy the conditions in
Theorems \ref{thmfila} and \ref{thmchara}.

We show that the conditions (1)
are satisfied.
Assume that the residue field
$F$ is perfect.
Then, we have
$G^n=G^{n-1}_{\rm cl}$
and 
$G^n_{\rm Ma}
=G^{n-1}_{\rm Ka}$.
Since $G^{n-1}_{\rm Ka}
=G^{n-1}_{\rm cl}$ in this case
by \cite[Th\'eor\`eme 9.11 (i)]{aml},
the condition (1)
in Theorem \ref{thmfila} is satisfied.

We identify 
$\Omega^1_{{\cal O}_K}
\otimes_{{\cal O}_K} F$
with 
${\mathfrak m}_K/
{\mathfrak m}_K^2$
by
$d\colon {\mathfrak m}_K/
{\mathfrak m}_K^2
\to 
\Omega^1_{{\cal O}_K}
\otimes_{{\cal O}_K} F$
and 
${\rm Hom}_F(
{\mathfrak m}_K^n/
{\mathfrak m}_K^{n+1},$
$\Omega^1_{{\cal O}_K}
\otimes_{{\cal O}_K}\bar F)$
with
${\rm Hom}_F(
{\mathfrak m}_K^{n-1}/
{\mathfrak m}_K^n,\bar F)$
by the induced isomorphism.
Then, the morphism
(\ref{eqgrMY})
is identified with
the morphism 
\begin{equation}
{\rm rsw}\colon
{\rm Hom}(G^{n-1}_{\rm cl}
/G^n_{\rm cl},{\mathbf F}_p)
\to
{\rm Hom}_F(
{\mathfrak m}_K^{n-1}/
{\mathfrak m}_K^n,F)
\label{eqgrK}
\end{equation}
defined in \cite[Corollary (5.2)]{K}.
Since the morphism
(\ref{eqgrK}) equals
(\ref{eqgrM}) by \cite[Th\'eor\`eme 9.11 (ii)]{aml},
the condition (1)
in Theorem \ref{thmchara} is satisfied.

We show that the conditions (2)
are satisfied.
For an extension 
$K'$ of henselian discrete valuation field
of ramification index $1$,
the diagram
\begin{equation}
\begin{CD}
{\rm Hom}(G^n_{\rm Ma}/
G^{n+1}_{\rm Ma},{\mathbf F}_p)
@>{{\rm rsw}'}>>
{\rm Hom}_F(
{\mathfrak m}_K^n/
{\mathfrak m}_K^{n+1},
\Omega^1_{{\cal O}_K}
\otimes_{{\cal O}_K}\bar F)\\
@VVV@VVV\\
{\rm Hom}(G^{\prime n}_{\rm Ma}
/G^{\prime n+1}_{\rm Ma},
{\mathbf F}_p)
@>{{\rm rsw}'}>>
{\rm Hom}_{F'}(
{\mathfrak m}_{K'}^n/
{\mathfrak m}_{K'}^{n+1},
\Omega^1_{{\cal O}_{K'}}
\otimes_{{\cal O}_{K'}}\bar F')
\end{CD}
\label{eqgrMM}
\end{equation}
is commutative.
Hence the condition (2) in Theorem \ref{thmchara}
is satisfied.

If $K'$ is tangentially dominant
over $K$,
then the morphism
$\Omega^1_{{\cal O}_K}
\otimes_{{\cal O}_K}\bar F
\to 
\Omega^1_{{\cal O}_{K'}}
\otimes_{{\cal O}_{K'}}\bar F'$
is an injection.
Hence by the commutative
diagram (\ref{eqgrMM}),
the morphism
$G^{\prime n}/G^{\prime n+1}
\to G^n/G^{n+1}$
is a surjection.
By the descending induction 
on $n$,
the condition (2) in Theorem \ref{thmfila}
is satisfied.
\qed
\end{proof}

\myaddress{School of Mathematical Sciences, University of Tokyo, Tokyo 153-8914, Japan}
\myemail{t-saito@ms.u-tokyo.ac.jp}


\begin{thebibliography}{99}

\bibitem{AS}
A.\ Abbes, T.\ Saito,
{\em Ramification of local fields with imperfect residue fields},
Amer.\ J.\ of Math., 124.5 (2002), 879-920.


\bibitem{aml}
A.\ Abbes, T.\ Saito,
{\em Analyse micro-locale $\ell$-adique en caract\'eristique $p>0$: le cas d'un trait},
Publ.\ RIMS 45 (2009), no.\ 1, 25-74.


\bibitem{Bo}
J.\ Borger,
{\em Conductors and the moduli of residual perfection},
Math.\ Ann.,
329 No.\ 1, (2004), 1-30.

\bibitem{K}
K.\ Kato,
{\em  Swan conductors for characters of degree 
one in the imperfect residue field case},  
Contemporary Math. {\bf 83}
(1989), 101--131. 

\bibitem{M}
S.\ Matsuda,
{\em On the Swan conductor in positive characteristic}, 
Amer.\ J.\ Math.\ 119 (1997), no.\ 4, 705--739.

\bibitem{TJM}
T.\ Saito,
{\em Ramification groups of coverings and valuations},
Tunisian J.\ of Math.,
Vol.\ 1, No.\ 3, 373-426, 2019.


\bibitem{red}
T.\ Saito,
{\em Graded quotients of ramification groups
of local fields with imperfect residue fields},
{\tt arxiv:2004.03768}.


\bibitem{CL}
J-P.~Serre,
{\sc Corps Locaux,}
Hermann, Paris, 1968.

\bibitem{Y}
Y.\ Yatagawa,
{\em Equality of two non-logarithmic ramification filtrations of abelianized Galois group in positive characteristic}, 
Doc.\ Math.\ 22 (2017), 917--952.
\end{thebibliography}
\end{document}